\documentclass[12pt]{article}

\usepackage{amsmath}
\usepackage{mathrsfs}
\usepackage{amsfonts}
\usepackage{graphicx}
\usepackage{subfig}
\usepackage{amssymb}

\begin{document}

\tableofcontents

\title{Linear Multistep Numerical Methods for Ordinary Differential Equations}
\author{Nikesh S. Dattani}

\maketitle

\begin{center}
Department of Applied Mathematics, University of Waterloo \\ Waterloo, Ontario~
 N2L 3G1, Canada
\end{center}

\begin{abstract}
A review of the most popular Linear Multistep (LM) Methods for solving Ordinary Differential Equations numerically is presented.  These methods are first derived from first principles, and are discussed in terms of their order, consistency, and various types of stability.  Particular varieties of stability that may not be familiar, are briefly defined first.  The methods that are included are the Adams-Bashforth Methods, Adams-Moulton Methods, and Backwards Differentiation Formulas.  Advantages and disadvantages of these methods are also described.  Not much prior knowledge of numerical methods or ordinary differential equations is required, although knowledge of basic topics from calculus is assumed.  

\end{abstract}

\section{Linear Multistep Methods}

\qquad As opposed to one-step methods, which only utilize one previous value of the numerical solution to approximate the subsequent value, multistep methods approximate numerical values of the solution by referring to \textit{more than one} previous value.  Accordingly, multistep methods may often achieve greater accuracy than one-step methods that use the same number of function evaluations, since they utilize \textit{more information} about the known portion of the solution than one-step methods do.

A special category of multistep methods are the \textit{linear} multi-step methods, where the numerical solution to the ODE at a specific location is expressed as a \textit{linear combination} of the numerical solution's values and the function's values at previous points.  For the standard system of ODEs, $\mbox{\boldmath $y'$} = \mbox{\boldmath $f$} (t,\mbox{\boldmath $y$)} $, a linear multistep method with k-steps would have the form:

\begin{equation}
\mbox{\boldmath $y_n$} = -\sum_{j=1}^k \alpha _j \mbox{\boldmath $y_{n-j}$} + h \sum_{j=0}^k \beta _j \mbox{\boldmath $f_{n-j}$},
\label{eq:generalLM}
\end{equation}

\noindent where $\alpha_j, \beta_j$ are constants, $\mbox{\boldmath $y_n$}$ is the numerical solution at $t = t_n$, and $\mbox{\boldmath $f_n$} = \mbox{\boldmath $f$}(t_n, \mbox{\boldmath $y_n$})$.   For the rest of this discussion, we will make the assumption that $f$ is differentiable as many times as needed, and we will consider the scalar ODE $y' = f(t,y)$ for simplicity in notation.  The generalization to systems of ODEs is straightforward.

It is important to note that in the above expression, all of the previous integration steps are assumed to be equally spaced, although it is possible to generalize these schemes to have variable step-sizes.  Also, note that if $\beta_0$ = 0 , the scheme is explicit (because it does not depend on $f_n$), and otherwise the scheme is implicit. We are now ready to examine some of the most popular families of linear multistep methods.
       
\subsection{The Adams Family}

\subsubsection{The Adams-Bashforth [AB] Methods (Explicit Adams Methods)}

\qquad  The most widely utilized linear multistep methods used for nonstiff problems are the Adams-Bashforth methods, which are members of the Adams Family that are explicit.  The spirit of the Adams-Bashforth technique is rooted in the Stone-Weierstrass Theorem.

\newtheorem{theorem}{Theorem}


\begin{theorem}
- Stone-Weierstrass Theorem - \emph{Let $f(t): \mathbb{R} \rightarrow \mathbb{C}$ be continuous on $ t \in [a,b]$.  For all $\epsilon > 0, \exists$ a polynomial $\phi(t) \ni ||f(t) - \phi(t)|| < \epsilon.$ }
\end{theorem}

\noindent In other words, any continuous function can be approximated to an arbitrary accuracy by a polynomial; generally, the more demanding the accuracy of the approximation, the higher the order needed of such a polynomial.  

With the Stone-Weierstrass Theorem in mind, we start with the ODE in question: $y' = f(t,y)$, and we integrate both sides to obtain:

\begin{equation}
y(t_n) = y(t_{n-1}) + \int_{t_{n-1}}^{t_n} f(t,y(t))dt.
\end{equation}


If we could integrate $f(t,y(t))$ analytically, we (likely) would not need to resort to numerical methods to determine the solution to the ODE.  If we \textit{cannot} integrate $f(t,y(t))$ analytically, according to the Stone-Weierstrass Theorem above, we can \textit{approximate} it with arbitrary accuracy by a polynomial $\phi(t)$, and since all polynomials can be integrated analytically, we have an obtainable, fair approximation of the solution to the ODE:

\begin{equation}
\label{eq:ABapprox}
y(t_n) \approx y(t_{n-1}) + \int_{t_{n-1}}^{t_n} \phi (t)dt.
\end{equation}

	Now to ensure that our approximation is reasonable, we insist that $\phi(t_{n-i}) = f(t_{n-i})$ for a reasonable number of integer values $i$.  For example, setting $\phi(t_{n-1})$ to be the constant $f(t_{n-1})$ will result in the scheme:
	
\begin{equation}
y(t_n) \approx y(t_{n-1}) + h f_{n-1}
\end{equation}	
	
Note that this scheme, which is also known as the 1-step Adams-Bashforth Method, is simply the classic Forward Euler (FE) method.  In terms of equation \eqref{eq:generalLM}, this scheme has $\alpha _1 = -1 , \beta _0 = 0 , \beta _1 = 1 $ and $ \beta _j, \alpha _j = 0$ for $j>1$.  

Let us now construct the 2-step Adams-Bashforth scheme.  We first need an interpolation polynomial $\phi (t)$ such that $\phi(t_{n-i}) = f(t_{n-1}) $ for $i = 1,2$.  The desired linear function is displayed below:

\begin{equation}
f(t,y) \approx \phi (t)= f(t_{n-2}) + \frac{f(t_{n-1}) - f(t_{n-2})}{t_{n-1} - t_{n-2}} (t - t_{n-2})
\end{equation}

\noindent Together with \eqref{eq:ABapprox}, we have that

\begin{equation}
y(t_n)  \approx  y(t_{n-1}) + \left[ f(t_{n-2}) t + \frac{f(t_{n-1}) - f(t_{n-2})}{t_{n-1} - t_{n-2}} \frac{(t - t_{n-2})^2}{2} \right]_{t_{n-1}}^{t_n} 
\end{equation}

\begin{equation}
y(t_n)  \approx  y(t_{n-1}) + h \left( \frac{3}{2} f(t_{n-1}) - \frac{1}{2}  f(t_{n-2}) \right)
\end{equation}

\noindent So in accordance with the line above, we can define our 2-step Adams-Bashforth scheme to be:

\begin{equation}
y_n = y_{n-1} + h ( \frac{3}{2} f_{n-1} - \frac{1}{2}f_{n-2}).
\end{equation}

\noindent This can also be expressed in terms of \eqref{eq:generalLM} with $\alpha _1 = -1 , \beta _0 = 0 , \beta _1 = \frac{3}{2} , \beta _2 = -\frac{1}{2} $ and $ \alpha _j = 0 , \beta _{j+1}$ for all $j > 1$.

Continuing in this manner, we can construct $k$-step Adams-Bashforth methods by interpolating $f$ through $k$ previous points: $t = t_{n-1} , t_{n-2} , \ldots t_{n-k}$.  Such a scheme could be derived by constructing a degree $\leq k-1$ polynomial $\phi (t)$ such that $\phi(t_{n-i}) = f(t_{n-i})$ for $i = 1, 2, \ldots k$, and integrating it as in \eqref{eq:ABapprox}, then replacing $y(t_n), y(t_{n-1})$ and $f(t_{n-i})$ with $y_n, y_{n-1} $ and $ f_{n-i}$ respectively.  As shown below, the resultant $k$-step Adams-Bashforth method can be expressed in the form of \eqref{eq:generalLM} with $\alpha _1 = -1, \beta _0 = 0 , \beta _j $ defined as displayed below for $1 \leq j \leq k$, and $\alpha _j = 0 , \beta _{j+k-1} = 0 $ for $j>1$:

\begin{equation}
y_n = y_{n-1} + h \sum_{j = 1}^{k} \beta _j f_{n-j} , 
\end{equation}
where
\begin{equation}
\beta _j = (-1)^{j-1} \sum_{i = j-1}^{k-1}\binom{i}{j-1}(-1)^i\int_{0}^{1}\binom{-s}{i}ds.
\end{equation}
 
It is important to mention that for such schemes, $k$ starting values must be given.  If only the initial condition is provided, the other $k-1$ points can be determined by a different scheme (for example, a Runge-Kutta method of the same order).  Adams-Bashforth methods also tend to have small regions of absolute stability (to be discussed later), and this inspired the construction of implicit Adams methods (called Adams-Moulton methods) which are the topic of the following discussion.   
 
\subsubsection{The Adams-Moulton [AM] Methods (Implicit Adams Methods)}

\qquad The difference between Adams-Moulton and Adams-Bashforth methods is that Adams-Moulton methods use an interpolating polynomial of degree $\leq k$ rather than $\leq k-1$ , and it includes $f$ at the unknown value $t_n$ as well.  A $k$-step Adams-Moulton scheme can be expressed in the form of \eqref{eq:generalLM} as follows:

\begin{equation}
y_n = y_{n-1} + h \sum_{j=0}^{k} \beta _j f_{n-j}.
\end{equation}

It is apparent that when $k=1$ and $\beta _1 = 0$ we have the classic Backward Euler (BE) method.  Likewise, if $k=1$ and $\beta _1 \neq 0$ we have the implicit trapezoidal method.

Adams-Moulton methods have smaller error constants, use less steps, and have larger stability regions than their Adams-Bashforth counterparts (of the same order).  However, AM methods using more than one step tend to have smaller regions of absolute stability than other implicit methods such as Runge-Kutta methods (in fact, they tend to be bounded, which often defeats the purpose of using an implicit scheme).  
Adams-Moulton methods have smaller error constants, use less steps, and have larger stability regions than their Adams-Bashforth counterparts (of the same order).  However, AM methods using more than one step tend to have smaller regions of absolute stability than other implicit methods such as Runge-Kutta methods (they tend to be bounded, which often defeats the purpose of using an implicit scheme).  
Adams-Moulton methods have smaller error constants, use less steps, and have larger stability regions than their Adams-Bashforth counterparts (of the same order).  However, AM methods using more than one step tend to have smaller regions of absolute stability than other implicit methods such as Runge-Kutta methods (they tend to be bounded, which often defeats the purpose of using an implicit scheme).  

  An alternative family of implicit linear multistep methods is the family of \textit{Backwards Differentiation Formulas}, which are the topic of the next section.  Such schemes are in fact the most popular methods for stiff problems.

\section{Backwards Differentiation Formulas [BDFs]}

\qquad	 In contrast to the linear multistep schemes in the Adams Family, who are derived by \textit{integrating} an interpolating polynomial $\phi(t)$ that approximates $f$, the BDF methods are derived by \textit{differentiating} an interpolating polynomial $\varphi(t)$ that approximates $y$ (one such that $\varphi(t_{n-i}) = y(t_{n-i})$ for $i = 0,1,2,\ldots k)$, and setting the derivative at $t_n$ to be equal to $f(t_n,y_n)$.  
	
	For example, the 1-step BDF method is derived as follows.  We first construct the interpolating polynomial $\varphi(t)$ that approximates $y$, with $\varphi(t_{n-i}) = y(t_{n-i})$ for $i = 0,1$.
	
\begin{equation}
y(t) \approx \varphi (t) = y(t_n) + (t - t_n) \frac{y(t_n) - y(t_{n-1})}{t_n - t_{n-1}}
\end{equation}
	
\noindent Upon differentiation, we get:

\begin{equation}
\label{eq:BDFdiff}
y'(t) = f(t,y) \approx \varphi '(t) = \frac{y(t_n) - y(t_{n-1})}{t_n - t_{n-1}}.
\end{equation}

\noindent We can then use the approximation in \eqref{eq:BDFdiff} as inspiration to construct our 1-step BDF method, by setting $\varphi '(t_n) = f(t_n,y_n)$:

\begin{equation}
\frac{y(t_n) - y(t_{n-1})}{h} = f(t_n,y_n).
\end{equation}

\noindent This is in fact, not surprisingly, equivalent to the Backward Euler method when rearranged.  

Similarly, we can construct a $k$-step BDF by generating the $k$-degree interpolating polynomial: 

\begin{equation}
\mbox{\scriptsize $y(t) \approx \varphi(t) = y_n + \frac{1}{h}(t-t_n)\nabla y_n + \frac{1}{2h^2}(t-t_n)(t-t_{n-1})\nabla ^2 y_n + \ldots + \frac{1}{h^k k!}(t-t_n)\ldots (t-t_{n-k+1})\nabla ^k, $}
\end{equation}

\noindent where $\nabla ^i$ is the backward difference operator:

\begin{equation}
\nabla ^0 y_n = y_n
\end{equation}

\begin{equation}
\nabla ^i y_n = \nabla ^{i-1} y_n - \nabla ^{i-1} y_{n-1}.
\end{equation}

\noindent Then upon differentiating, and setting $\varphi '(t_n) = f(t_n,y_n)$ we get:

\begin{equation}
\sum_{i=1}^{k} \frac{1}{i} \nabla ^i y_n = h f(t_n,y_n),
\end{equation}

\noindent which can be transformed to match the general expression of \eqref{eq:generalLM}, with $\beta _j = 0 $ for $j>0$ (note that this makes them \textit{implicit} schemes):

\begin{equation}
y_n = -\sum_{i=1}{k} \alpha _i y_{n-i} + h \beta _0 f_n.
\end{equation}

It is noteworthy that the backward differences $\nabla ^i y$ of $y$ approximate the true derivatives of $y$ (i.e. $\nabla ^i y \approx h^k y^{(i)}$).

What makes BDF methods powerful is their unique and convenient region of absolute stability giving rise to its L-stability (to be discussed later), which is of particular importance for stiff problems.  It is for this same reason that BDF methods with $k>6$ are not used, as these methods have a region of absolute stability that crosses the negative real axis, thereby disqualifying them from being classified as A-stable.  

\section{Order and consistency of Linear Multistep Methods}

\qquad  Investigating convergence of linear multistep methods is quite different from that of non-linear one-step methods (such as the Runge-Kutta) methods.  For Runge-Kutta methods, 0-stability is automatic, and investigating the order can be cumbersome.  For linear multistep methods, 0-stability is not necessarily automatic, and needs to be confirmed for each scheme.  Contrarily, and again unlike the Runge-Kutta methods, investigating the order of linear multistep methods is rather straightforward.  We will start this section by investigating the order of linear multistep methods.

\subsection{Order}

To begin, let us define the linear operator $\mathscr{L}_h [y(t)]$, where $y(t)$ is an arbitrarily continuously differentiable funciton on [0,b]:

\begin{equation}
\mathscr{L}_h [y(t)] = \sum_{j=0}^{k}[\alpha_j y(t-jh) - h \beta _j y'(t-jh)]
\end{equation} 

\noindent This expression is based on Eq. \eqref{eq:generalLM}.  Recalling that $y' = f(t,y(t))$, we can write the above expression in the following way:

\begin{equation}
\mathscr{L}_h [y(t)] = \sum_{j=0}^{k}[\alpha_j y(t-jh) - h \beta _j f(t-jh,y(t-jh))]
\end{equation}

\noindent which becomes, after expanding $y(t-jh)$ and $f(t-jh,y(t-jh))$ in a Taylor series about $t$ and simplifying:

\begin{equation}
\label{eq:discretizationOperatorTaylor}
\mathscr{L}_h [y(t)] = C_0 y(t) + C_1 h y'(t) + \ldots + C_q h^q h^{q}(t) + \ldots,
\end{equation}

\noindent where,

\begin{equation}
C_0 = \sum_{j = 0}^{k}\alpha _j, \textrm{and} 
\end{equation}

\begin{equation}
C_i = (-1)^i \left[ \frac{1}{i!} \sum_{j=1}^{k} j^i \alpha _j + \frac{1}{(i-1)!} \sum_{j=0}^{k}j^{i-1} \beta _j \right] , i = 1,2,3, \ldots
\end{equation}

Now, the order of the method is $p$ if the local truncation (or discretization) error is $d_n = O(h^p)$, which is given by:

\begin{equation}
d_n = \frac{\mathscr{L}_h [y(t_n)]}{h},
\end{equation}

\noindent where $y(t_n)$ is the exact solution at $t = t_n$.  So using the version of the expression given by Eq. \eqref{eq:discretizationOperatorTaylor}, we get that the method is order $p$ if:


\begin{equation}
C_0 = C_1 = \ldots = C_p = 0 , C_{p+1} \neq 0.
\end{equation}

\noindent Combining this result with \eqref{eq:discretizationOperatorTaylor}, we get that 

\begin{equation}
d_n = C_{p+1}h^py^{(p+1)}(t_n) + O(h^{p+1}),
\end{equation} 

\noindent where $C_{p+1}$ is the \textit{error constant} of the scheme.  

It can then be shown that Adams-Bashforth and BDF methods are of order $k$ (where $k$ is the number of steps), while Adams-Moulton methods are of order $k+1$ (with the exception of the case where the scheme is completed in a single step with $\beta _1 = 0$, as in Backward Euler, which is order $k=1$).  

Calculating each value of $C_q$ in \eqref{eq:discretizationOperatorTaylor} can be tedious though.  The easiest way to check for consistency ($p \geq 1$) is to use the fact that a method is consistent \textit{if and only if}

\begin{equation}
\sum_{j=0}^{k} \alpha _j = 0 \textrm{ and } \sum_{j=1}^{k}j \alpha _j + \sum_{j=0}^{k} \beta _j = 0.
\end{equation} 

This result can also be demonstrated in terms of the \textit{characteristic polynomials} of the recurrence relations arising from the expression for the numerical scheme:

\begin{equation}
\rho(\xi) = \sum_{j=0}^{k}\alpha _j \xi ^{k-j}, \alpha _0 \equiv 1,
\end{equation}

\begin{equation}
\sigma(\xi) = \sum_{j=0}^{k}\beta _j \xi ^{k-j}.
\end{equation}

In particular, the scheme is consistent \textit{if and only if} $\rho (1) = 0$, $\rho '(1) = \sigma(1)$.

\section{Stability of Linear Multistep Methods}

\subsection{0-Stability}

\qquad Here we are looking at Eq. \eqref{eq:generalLM} in the limit as $h \rightarrow 0$, rendering this equation into the form:

\begin{equation}
\label{eq:characEq}
\alpha _k y_{n-k} + \alpha _{k-1} y_{n-k+1} + \ldots \alpha _0 y_n = 0.
\end{equation}

\noindent This recurrence relation has the characteristic polynomial $\rho (\xi)$ defined earlier.  Due to consistency (see previous section), we know that $\xi = 1$ is a root of Eq. \eqref{eq:characEq}.  If the rest of the roots of the equation are distinct, then we have a solution of the form:

\begin{equation}
y_n = \sum_{i=1}^{k-1}c_i + \xi _{i}^n + c_{k} (1)^n.
\end{equation}

\noindent If $\xi _1 = \xi _2$ is a double root, then the solution will have the form:

\begin{equation}
y_n = \sum_{i=3}^{k-1}c_i + \xi _{i}^n + c_{k} (1)^n + c_{1} \xi _1^n + c_2 n \xi _2^n.
\end{equation}

\noindent Similarly, if $\xi _1 = \xi _2 = \xi _3$ is a triple root, then the solution will have the form:

\begin{equation}
y_n = \sum_{i=4}^{k-1}c_i + \xi _{i}^n + c_{k} (1)^n + c_{1} \xi _1^n + c_2 n \xi _2^n + c_3 n (n-1) \xi _3^n.
\end{equation}

\noindent We can see that if $|\xi _i| > 1$ , then our solution will diverge as $n$ gets large.  Likewise, if $|\xi _i| = 1$ is not a simple root, we will again have divergence.  Therefore, a linear multistep method is 0-stable \textit{if and only if} all roots of the equation $\rho (\xi) = 0$ satisfy $|\xi _i| \leq 1$ , where if $|\xi _i| = 1$, then $\xi _i$ is a simple root, for $1 \leq i \leq k$.  Now by a theorem sometimes referred to as the Dahlquist Theorem, if this root condition is satisfied, \textit{and} the method is accurate to order $p$, and the initial values are accurate to order $p$, then the method is convergent to the order $p$.

We can further specify the \textit{strength} of the stability of the scheme by defining \textit{strongly stable} as meaning all roots of $\rho(\xi)=0$ have the property $|\xi| < 1$ with the exception of the one root $\xi _i$ which equals 1.  A scheme can then be defined as weakly stable if it is 0-stable, but \textit{not} strongly stable.

\subsection{Absolute Stability, A-stability and L-stability}

\qquad By applying Eq. \eqref{eq:generalLM} to the test equation $y' = \lambda y$ and letting $y_n = \xi ^n$, we find that $\xi$ must satisfy $\rho(\xi) - h\lambda \sigma(\xi)$.  Now to address absolute stability, we define the stability polynomial as $\zeta = \rho(\xi) - h \lambda \sigma(\xi)$ and note that the absolute stability region is the region of values of $\lambda h$ such that $|y_n|$ does not grow with increasing values of $n$.  For this condition to be met, all roots $\xi _i$ of $\zeta = 0$ must satisfy $|\xi_i| \leq 1$.  


As for A-stability, recall that a numerical scheme is A-stable if its absolute stability region contains the entire left half of the complex plane (i.e. it contains Re($\lambda h <$ 0)).  It turns out that explicit linear multistep methods cannot be A-stable, and that A-stable linear multistep methods with order greater than 2 do not exist.  The most accurate LM method (method with smallest error constant) is the second-order implicit trapezoidal method (with error constant $C_3 = \frac{1}{12}$).  We conclude that A-stability is rare in linear multi-step methods.

After applying our numerical scheme to the test equation $y' = \lambda y$, we can rearrange our expression to obtain an expression of the form $y_n = R(z)y_{n-1}, $ where $z=\lambda h$.  If $R(z) \rightarrow 0$ as $Re(z) \rightarrow -\infty$, we say that the scheme is L-stable (or has \textit{stiff decay}.  This explains why the Backward Euler method works better than the Trapezoidal Rule for some problems, even though the Trapezoidal Rule is of higher order!

\newpage





\end{document}